\documentclass[a4paper, 10pt]{article}
 
\usepackage{authblk}
\usepackage[latin1]{inputenc}
\usepackage[T1]{fontenc}
\usepackage{amsfonts}
\usepackage{graphicx}
\usepackage{ulem}
\usepackage{lmodern}
\usepackage[top=3cm, bottom=3cm, left=4cm, right=4cm]{geometry}
\usepackage{amsthm}
\usepackage{amsmath}
\usepackage{amssymb}
\usepackage{mathrsfs}
\usepackage{enumerate}
\usepackage{sectsty}
\newtheorem*{Definition}{Definition}
\newtheorem{Lemma}{Lemma}
\newtheorem{Theorem}{Theorem} 
\newtheorem{Corollary}{Corollary}  
\usepackage{fancyhdr}
\allsectionsfont{\centering}
\begin{document}
\title{A THEOREM FOR DISTINCT ZEROS OF $L$-FUNCTIONS}
\author{Quentin Gazda\footnote{quentin.gazda@ens-cachan.fr}}
\affil{\'Ecole Normale Sup\'erieure}
\affil{94230 Cachan}
\date{\today} 
\maketitle
\begin{abstract} In this paper, we establish a simple criterion for two $L$-functions $L_1$ and $L_2$ satisfying a functional equation (and some natural assumptions) to have infinitely many distinct zeros. Some related questions have already been answered in the particular case of Automorphic forms using so-called Converse Theorems (see \cite{Rag}). Deeper results can also be stated for elements of the Selberg class (see \cite{Bom}). However, we shall give here a general answer that do not use any advanced topics in analytic number theory. Therefore, this paper should be accessible to anyone who has some basic notions in measure-theory and advanced complex analysis.
\end{abstract}
Throughout this paper, we distinguish two kinds of Dirichlet series, as it is usually done. A meromorphic function $f$ on a domain $\Omega \subseteq \textbf{C}$ is represented by a General Dirichlet Series (GDS) if $f$ can be written in the form
\begin{equation}
f(s)=\sum_{n=1}^{\infty}{a_n e^{-\lambda_n s}} \nonumber
\end{equation}
($(a_n)$ is a sequence of complex numbers and $(\lambda_n)$ is a sequence of real and non-negative numbers strictly increasing to infinity) when $\operatorname{Re}(s)$ is large enough. In the same framework, $f$ is represented by an Ordinary Dirichlet Series (ODS) if $f$ can be written in the form
\begin{equation}
f(s)=\sum_{n=1}^{\infty}{\frac{a_n}{n^s}} \nonumber
\end{equation}
when $\operatorname{Re}(s)$ is large enough. In the event that the GDS is absolutely convergent (when $\operatorname{Re}(s)$ is large), we shall say that this is an absolutely convergent GDS. In this case, we note $\sigma_a(f)$ the abscissa of absolute convergence of the series (this number depends only on $f$ and not on the form of the Dirichlet series as we will see, Corollary \ref{uniq}). We recall that for ODS, the finite abscissa of convergence is equivalent to the finite abscissa of absolute convergence via the inequality $\sigma_c(f)\leq \sigma_a(f) \leq \sigma_c(f)+1$ (where $\sigma_c(f)$ is the abscissa of convergence of the ODS).

The purpose of this paper is to provide a relationship between two Dirichlet series which may be extended to a meromorphic function of finite growth on the whole complex plane, satisfying a functional equation and whose sets of zeros (counted with multiplicity) differ from one another only by a finite number of elements. More precisly, we shall prove
\begin{Theorem} \label{main}
Let $\gamma$ be a meromorphic function on $\textbf{C}$, $k$ a real. Suppose that $L_1$ and $L_2$ are two meromorphic functions on $\textbf{C}$ such that $L_i$ (for $i=1$ and $2$) satisfies the following conditions
\begin{enumerate}[(i)]
\item There exists a meromorphic function $L_i^*$ on $\textbf{C}$ and a positive real $N_i$ such that
\begin{equation}
L_i(k-s)=N_i^s\gamma(s)L_i^*(s) \quad (\forall s\in\textbf{C}), \nonumber
\end{equation}
\item there exists a polynomial $P_i$ such that $P_iL_i$ is holomorphic and of finite order in $\textbf{C}$. 
\end{enumerate}
Assume further that $L_1/L_2$ has only finitely many poles and that $L_1/L_2$ and $L_1^*/L_2^*$ are represented by ODS.

Then there exist an integer $N$ and complex numbers $(a_d)_{d|N}$  such that
\begin{equation}
L_1(s)=\left(\sum_{u|N}{\frac{a_u}{u^s}}\right)L_2(s) \quad (\forall s\in \textbf{C}) \nonumber
\end{equation}
\end{Theorem}

The proof will infer that $N=N_1/N_2$ when $N_1/N_2$ is actually a positive integer. The latter result can be seen as a generalization of a famous result of Hamburger (see \cite{Ham}). It has the following well-known corollaries (we left the proof for the reader)
\begin{Corollary}
Let $\chi$ and $\varphi$ be two distinct primitive Dirichlet characters. Then $L(\chi,s)/L(\varphi,s)$ has infinitely many poles. 
\end{Corollary}

\begin{Corollary}
Let $f$ and $g$ be two linearly independent Hecke forms of integer weight $k$ for $\Gamma_0(N)$. Then $L(f,s)/L(g,s)$ has infinitely many poles.
\end{Corollary}
\section{FUNCTIONS REPRESENTED BY GDS}
We recall that for a measure $\mu$ on the measurable space $([0,\infty),\mathcal{B}[0,\infty))$ (as usual, $\mathcal{B}[0,\infty)$ denotes the Borel sets of the topological space $[0,\infty)$), the Laplace transform $\mathcal{L}\mu$ of $\mu$ is defined by :
\begin{equation}
\mathcal{L}\mu(s):=\int_{0}^{\infty}{e^{-st}\mu(dt)}. \nonumber
\end{equation}
When the image measure $|\mu|$ has finite mass, $\mu$ will be refered to as absolutely convergent : in this case, the Laplace transform of $\mu$ converges for $\operatorname{Re}(s)\geq 0$. 

The Laplace transform will be useful for us through the following example : let $\nu$ be the measure on $([0,\infty),\mathcal{B}[0,\infty))$ defined by
\begin{equation}
\nu=\sum_{n=1}^{\infty}{a_n\delta_{\lambda_n}}. \nonumber
\end{equation} 
It is straightforward to see that the corresponding GDS is the Laplace transform of $\nu$ (remark that the measure is absolutely convergent if and only if $\sigma_a(\mathcal{L}\nu)$ is non-positive). Therefore, via a uniqueness theorem for Laplace transform, we should be able to derive an equivalent result for GDS.

\begin{Theorem} \label{lap}
Let $\mu_1$ and $\mu_2$ be two absolutely convergent measures on $([0,\infty),\mathcal{B}[0,\infty))$ such that $\mathcal{L}\mu_1(s)=\mathcal{L}\mu_2(s)$ for all $s\geq 0$. Then $\mu_1=\mu_2$.
\end{Theorem}

Note that the assumption : $\mathcal{L}\mu_1(s)=\mathcal{L}\mu_2(s)$ for all $s\geq 0$ can be weaken in : $\mathcal{L}\mu_1(s)=\mathcal{L}\mu_2(s)$ for all $s\geq \sigma$ (for a $\sigma \in \textbf{R}$),
if we assume that the measure $(t\mapsto e^{-\sigma t})\mu_i$ is absolutlely convergent (for $i=1$ and $2$).

\paragraph{Proof :} Setting $\mu:=\mu_1-\mu_2$ (also absolutely convergent), the assumption is equivalent to
\begin{equation}
\int_{0}^{\infty}{e^{-st}\mu(dt)}=0 \quad (\forall s\geq 0). \nonumber
\end{equation}
By the Stone-Weierstrass theorem, for any continous function $f$ on $[0,1]$, we can find a sequence $(P_n)$ of polynomials such that
\begin{equation}
\|f-P_n\|_{\infty}:=\operatorname{sup}_{t\in [0,1]}|f(t)-P_n(t)| \underset{n\to \infty}{\longrightarrow} 0. \nonumber
\end{equation}
Hence,
\begin{equation}
\left|\int_{0}^{\infty}{f(e^{-t})\mu(dt)} \right|=\left|\int_{0}^{\infty}{f(e^{-t})-P_n(e^{-t})\mu(dt)} \right|\leq |\mu|([0,\infty))\|f-P_n\|_{\infty} \nonumber
\end{equation}
and since the right-hand side is arbitrary small (for large values of $n$), we conclude that 
\begin{equation}
\int_{0}^{\infty}{f(e^{-t})\mu(dt)}=0 \nonumber
\end{equation}
for any continous function $f$ on $[0,1]$. Therefore, we fix $a<b$ in $(0,\infty)$ and we define $f_{\varepsilon}$ (for $\varepsilon>0$) by
\begin{equation}
f_{\varepsilon}(x):=\begin{cases} 
      \frac{1}{\varepsilon}(x-e^{-a}+\varepsilon) & e^{-a}-\varepsilon \leq x \leq e^{-a} \\
      1 & e^{-a}\leq x\leq e^{-b} \\
      -\frac{1}{\varepsilon}(x-e^{-b}-\varepsilon) & e^{-b} \leq x \leq e^{-b}+\varepsilon \\
			0 & \text{otherwise.}
   \end{cases} \nonumber
\end{equation}
If $\varepsilon$ is small enough, $f_{\varepsilon}$ is continous and defined on $[0,1]$. Thus,
\begin{equation}
|\mu([a,b])|=\left|\int_{0}^{\infty}{\textbf{1}_{[a,b]}(t)\mu(dt)}\right|\leq \left|\int_{0}^{\infty}{f_{\varepsilon}(e^{-t})\mu(dt)}\right|+2|\mu|([0,\infty))\varepsilon=2|\mu|([0,\infty))\varepsilon \nonumber
\end{equation}
and finally $\mu([a,b])=0$ for all $b>a\geq 0$ (the case $a=0$ can be done similary). To conclude, we define
\begin{equation}
\mathcal{C}:=\left\{[a,b] | a,b \in [0,\infty) \right\}\subset \mathcal{B}[0,\infty). \nonumber
\end{equation}
$\mathcal{C}$ is a $\pi$-system such that $\sigma(\mathcal{C})=\mathcal{B}[0,\infty)$ and on which $\mu$ and the trivial measure coincide. Finally, knowing that $\mu([0,\infty))=0$, a well-known result from measure-theory allows us to conclude that $\mu$ is the trivial measure. Therefore $\mu_1=\mu_2$. \\

This theorem has many consequences in the theory of GDS. It says in particular that a very few functions are represented by an absolutely convergent GDS.

\begin{Corollary}[Uniqueness Theorem for GDS] \label{uniq}
Let $f:[a,\infty)\rightarrow\textbf{C}$ for a certain $a<\infty$. Then $f$ can be represented by at most one absolutely convergent GDS. 
\end{Corollary}

\paragraph{Proof :} Suppose that we can write (for $s$ large enough)
\begin{equation}
f(s)=\sum_{n=1}^{\infty}{a_ne^{-\lambda_n s}}=\sum_{n=1}^{\infty}{b_ne^{-\mu_n s}} \nonumber
\end{equation}
where the two series are absolutely converging for a certain $\sigma\in \textbf{R}$. Without loss of generality, we may and shall assume that $\sigma=0$ (taking $a'_n:=a_ne^{-\sigma \lambda_n}$ and $b'_n:=b_ne^{-\sigma \mu_n}$). We introduce
\begin{equation}
\nu_1=\sum_{n=1}^{\infty}{a_n\delta_{\lambda_n}}, \quad \nu_2=\sum_{n=1}^{\infty}{b_n\delta_{\mu_n}}. \nonumber
\end{equation}
Those two measures are absolutely convergent (since $\sigma=0$) and $\mathcal{L}\nu_1(s)=\mathcal{L}\nu_2(s)$ for all $s\geq 0$. The previous theorem implies that $\nu_1=\nu_2$. Therefore, by measuring the singletons $(\left\{\lambda_n\right\})$, $(\left\{\mu_n\right\})$ under this equality, one obtains that the two series are the same.

\begin{Corollary} \label{deux}
A non-constant rational fraction is not represented by an absolutely convergent GDS.
\end{Corollary}

\paragraph{Proof :} We assume the opposite, i.e. there exists $F$ a non-constant rational function represented by an absolutely convergent GDS. For $s\geq \sigma_a(F)$,
\begin{equation}
F(s)=\frac{a_n s^n+...+a_0}{b_m s^m+...+b_0}=\sum_{n=1}^{\infty}{a_n e^{-\lambda_n s}}. \nonumber
\end{equation}
From the partial fraction principle, we know that we can write $F$ as
\begin{equation}
F(s)=\frac{a_n s^n+...+a_0}{b_m s^m+...+b_0}=P(s)+F_1(s)+F_2(s)+...+F_p(s) \nonumber
\end{equation}
where $P$ is a polynomial and where (for $i \in \left\{1,...,p\right\}$),
\begin{equation}
F_i(s)=\frac{c_{1,i}}{s-z_i}+\frac{c_{2,i}}{(s-z_i)^2}+...+\frac{c_{n_i,i}}{(s-z_i)^{n_i}} \nonumber
\end{equation}
for complex numbers $(z_i)$, $(c_{i,j})$ and positive integers $(n_i)$. Now, by the absolute convergence of the GDS, $\lim_{s\to \infty}F(s)$ is finite, we denote it $\ell$. $P$ must be constant equal to $\ell$ since $F_i$ (for $i \in \left\{1,...,p\right\}$) tends to zero at infinity. 

Let $f_i$ (for $i \in \left\{1,...,p\right\}$) be the function
\begin{equation}
f_i(t):=\sum_{r=1}^{n_i}{c_{r,i}\frac{t^{r-1}e^{z_i t}}{(r-1)!}} \nonumber
\end{equation}
(whose Laplace transform is $F_i$) and let $\lambda$ be the Lebesgue measure on $\textbf{R}$. In the same framework, Theorem \ref{lap} implies that the two measures
\begin{equation}
\ell \delta_0+\sum_{i=1}^p{f_i \lambda} \quad \text{and} \quad \sum_{n=1}^{\infty}{a_n \delta_{\lambda_n}} \nonumber
\end{equation}
are the same. We easily see that $f_i$ (for $i \in \left\{1,...,p\right\}$) must be zero. A contraction since $F$ is non-constant.
\section{THE SPACE $\mathcal{O}_N$}
For convenience, we introduce a definition : a meromorphic function $f$ is of \textbf{finite order} in a set $\Omega \subseteq \textbf{C}$ if a $\rho>0$ exists such that $f(s)\ll e^{|s|^{\rho}}$ in $\Omega$. It extends the usual definition of the finite order for entire functions (taking $\Omega=\textbf{C}$).

\begin{Definition}
For $N$ a positive real, let $\mathcal{O}_N$ (resp. $\mathcal{G}_N$) be the space of functions $f$ that satisfy :
\begin{enumerate}
\item $f$ is meromorphic on $\textbf{C}$ and represented by an ODS (resp. absolutely convergent GDS),
\item $\forall s\in \textbf{C}$, $f(-s)=N^s g(s)$ where $g$ is represented by an ODS (resp. absolutely convergent GDS),
\item $M>0$ exists such that $f$ is of finite order in the set
\begin{equation}
\Omega_M(f):=\left\{-1-\sigma_a(g)\leq \operatorname{Re}(s) \leq \sigma_a(f)+1 \right\}\cap \left\{|\operatorname{Im}(s)|\geq M\right\}. \nonumber
\end{equation} 
\end{enumerate}
\end{Definition}

It is clear that if $f$ belongs to $\mathcal{G}_N$ (resp. $\mathcal{O}_N$) then $g$ belongs to $\mathcal{G}_{N}$ (resp. $\mathcal{O}_{N}$). The goal of this section is to understand the structure of $\mathcal{O}_N$.
\begin{Lemma} \label{dc}
$f$ belongs to $\mathcal{G}_1$ if and only if $f$ is constant. In particular $\mathcal{G}_1=\mathcal{O}_1=\operatorname{Vect}\left\{ 1 \right\}$.
\end{Lemma}

\paragraph{Proof :} We split the complex plane into the following sets
\begin{equation}
\Omega_1:=\Omega_M(f) \nonumber
\end{equation}
\begin{equation}
\Omega_2:=\left\{\operatorname{Re}(s) \geq \sigma_a(f)+1\right\}  \nonumber
\end{equation}
\begin{equation}
\Omega_3:=\left\{\operatorname{Re}(s)\leq -1-\sigma_a(g)\right\}  \nonumber
\end{equation}
\begin{equation}
\Omega_4:=\left\{-1-\sigma_a(g)\leq\operatorname{Re}(s)\leq \sigma_a(f)+1\right\}\cap\left\{|\operatorname{Im}(s)|\leq M\right\}. \nonumber
\end{equation}
By the fact that $f(-s)=g(s)$ and the general theory of Dirichlet series, we know that $f$ is bounded on $\Omega_2$ and $\Omega_3$. In particular, $f$ is bounded on the lines $\left\{\operatorname{Re}(s)=\sigma_a(f)+1\right\}$ and $\left\{\operatorname{Re}(s)=-1-\sigma_a(g)\right\}$. We can apply the Phragm\`en-Lindel\"of principle to $\Omega_1$ to conclude that $f$ is bounded on this set. Thus, since $\Omega_1 \cup\Omega_2 \cup\Omega_3 \cup\Omega_4=\textbf{C}$, the only poles of $f$ are in $\Omega_4$ which is bounded. By the isolated zeros principle, $f$ has only finitely many poles. Let $\rho_1$, $\rho_2$,..., $\rho_m$ be those poles counted with multiplicity. To sum up, we can find $A$ and $B$ two positive constants such that
\begin{equation}
\left| \prod_{i=1}^m{(s-\rho_i)}f(s) \right|\leq A+B|s|^m \quad (\forall s\in \textbf{C}). \nonumber
\end{equation}
The general form of Liouville's theorem then implies that the left hand side is a polynomial in $s$. Hence, $f$ is a rational fraction in $s$. By Corollary \ref{deux}, $f$ must be constant.

\begin{Lemma}
If $N>0$ is not an integer, $\mathcal{O}_N=\left\{0 \right\}$.
\end{Lemma} 

\paragraph{Proof :} Clearly, the constant function equal to zero belongs to $\mathcal{O}_N$. We now establish the reverse inclusion. We take $f\in \mathcal{O}_N$ and write $\left\lfloor N \right\rfloor\neq N$ the entire part of $N$. We define
\begin{equation}
\tilde{f}:s \longmapsto N^{s/2}\left(f(s)-\sum_{u=1}^{\left\lfloor N \right\rfloor}{\frac{a_u}{u^s}}\right)=\sum_{n=1}^{\infty}{a_{\left\lfloor N \right\rfloor+n}e^{-s(\log(\sqrt{N})-\log(\left\lfloor N \right\rfloor+n))}}. \nonumber
\end{equation}
An easy computation and the above formula show that $\tilde{f}$ belongs to $\mathcal{G}_1$ and as a consequence is constant by Lemma \ref{dc}. Since $\log(\left\lfloor N \right\rfloor+1)> \log(\sqrt{N})$ (because $N>1$) and the absolute convergence of the series (for $\operatorname{Re}(s)$ large), $\tilde{f}(s)$ tends to zero as $s$ tends to infinity : this constant must be zero. It gives
\begin{equation}
f(s)=\sum_{u=1}^{\left\lfloor N \right\rfloor}{\frac{a_u}{u^s}} \quad (\text{for~$\operatorname{Re}(s)$~large}). \nonumber
\end{equation}
Since $N^{-s}f(-s)=g(s)$ ($\forall s\in \textbf{C}$), we have
\begin{equation}
g(s)=\sum_{u=1}^{\left\lfloor N \right\rfloor}{a_u\left(\frac{u}{N} \right)^s} \quad (\text{for~$\operatorname{Re}(s)$~large}) \nonumber
\end{equation}
and thus everywhere on $\textbf{C}$ by analytic continuation. But, by assumption, $g$ is represented by an ODS. From Corollary \ref{uniq}, since the right-hand side is absolutely convergent (as a finite sum), we obtain that this right-hand side must be an ODS. A contradiction since $N$ is not an integer. \\

The complete structure of $\mathcal{O}_N$ is given by the following theorem :
\begin{Theorem} \label{even}
For all $N>0$, the space $\mathcal{O}_N$ is a finite-dimensional $\textbf{C}$-vector space and when $N$ is an integer,
\begin{equation}
\mathcal{O}_N=\operatorname{Vect}_{u|N}\left\{s \longmapsto \frac{1}{u^s}\right\}. \nonumber
\end{equation}
In particular, $\operatorname{dim}(\mathcal{O}_N)=d(N)$ (where $d(N)$ denotes the number of positive divisors of $N$ if $N$ is an integer and $0$ otherwise). 
\end{Theorem}

\paragraph{Proof :} Without loss of generality, we assume that $N\geq 1$ is an integer (the other case is given by the previous lemma). It follows directly from the definition that $\mathcal{O}_N$ is a $\textbf{C}$-vector space. Since $u^s=\left(\frac{u^s}{N^s}\right)N^s$, for $u$ a positive divisor of $N$, we have
\begin{equation}
s \longmapsto \frac{1}{u^s} \in \mathcal{O}_N \nonumber
\end{equation}
and thus,
\begin{equation}
\operatorname{Vect}_{u|N}\left\{s \longmapsto \frac{1}{u^s}\right\} \subseteq \mathcal{O}_N \subseteq \mathcal{G}_N \nonumber
\end{equation}
Similary, we also have
\begin{equation}
\operatorname{Vect}_{u\in \left\{1,...,N\right\}}\left\{s \longmapsto \frac{1}{u^s}\right\} \subseteq \mathcal{G}_N \nonumber
\end{equation}
We turn on to the other inclusion. Let $f(s)=\sum_{n\geq 1}{a_n n^{-s}} \in \mathcal{O}_N \subseteq \mathcal{G}_N$ and $\tilde{f}$ be the function defined as
\begin{equation}
\tilde{f}:s \longmapsto N^{s/2}\left(f(s)-\sum_{u=1}^N{\frac{a_u}{u^s}}\right)=N^{s/2}\sum_{n>N}{\frac{a_u}{u^s}} \nonumber
\end{equation}
which is in $\mathcal{G}_N$ thanks to the vector-space structure. By construction, $\tilde{f}$ admits a GDS of the form
\begin{equation}
\tilde{f}(s)=\sum_{n=1}^{\infty}{b_n e^{-s \mu_n}}=\sum_{n=1}^{\infty}{a_{N+n}e^{-s(\log(\sqrt{N})-\log(N+n))}} \nonumber
\end{equation}
where $(\mu_n)$ is a strictly increasing sequence of positive numbers (as in the proof of the previous lemma). Consequently, $s \longmapsto \tilde{f}(s)$ lies in $\mathcal{G}_1$. By Lemma \ref{dc}, $\tilde{f}$ is constant on $\textbf{C}$ : this constant must be zero by construction, since $f$ has no constant term. We have :
\begin{equation}
f(s)=\sum_{u=1}^N{\frac{a_u}{u^s}} \nonumber
\end{equation}
Finally $f$ must belong to $\mathcal{O}_N$, so we should add the condition $u|N$ in the index of the sum (this can be seen by invoking Corollary \ref{uniq}, again). Thus
\begin{equation}
\operatorname{Vect}_{u|N}\left\{s \longmapsto \frac{1}{u^s}\right\} \supseteq \mathcal{O}_N \nonumber
\end{equation}
which proves the result.
\section{PROOF OF THEOREM \ref{main}}
It remains to prove a lemma.
\begin{Lemma} \label{order}
Let $f$ and $g$ be two entire functions of finite order such that $f/g$ is entire. Then $f/g$ is of finite order.
\end{Lemma}

\paragraph{Proof :}
Since $f$ and $g$ are of finite order, they both have finite genus (see \cite{Gar} for a general study of the order, the genus and Hadamard's factorization theorem). Thus, they both admit a Hadamard product :
\begin{equation}
f(s)=s^{n}e^{P(s)}\prod_{\rho \in \mathcal{Z}(f)}{\left(1- \frac{s}{\rho}\right)\exp\left({\frac{s}{\rho}+\frac{s^2}{2\rho^2}+...+\frac{s^{h}}{h\rho^{h}}}\right)} \nonumber
\end{equation}
\begin{equation}
g(s)=s^{m}e^{Q(s)}\prod_{\rho \in \mathcal{Z}(g)}{\left(1- \frac{s}{\rho}\right)\exp\left({\frac{s}{\rho}+\frac{s^2}{2\rho^2}+...+\frac{s^{r}}{r\rho^{r}}}\right)} \nonumber
\end{equation}
where we use the notation $\mathcal{Z}(f)$ for the multiset of zeros of $f$. Here, $P$ and $Q$ are two polynomials, and $h$ (resp. $r$) is the genus of $f$ (resp. $g$). By assumption $\mathcal{Z}(g)\subset \mathcal{Z}(f)$ (where the inclusion for multisets means that all elements of $\mathcal{Z}(g)$ are in $\mathcal{Z}(f)$ with a lower multiplicity). Consequently, the ratio of the two products over zeros of $f$ and $g$ can be simplified in
\begin{equation}
\exp\left(\sum_{i=r+1}^{h}{\frac{s^i}{i}\left(\sum_{\rho \in \mathcal{Z}(g)}{\frac{1}{\rho^i}}\right)}\right)\prod_{\rho \in \mathcal{Z}(f)\backslash \mathcal{Z}(g)}{\left(1- \frac{s}{\rho}\right)\exp\left({\frac{s}{\rho}+\frac{s^2}{2\rho^2}+...+\frac{s^{h}}{h\rho^{h}}}\right)} \nonumber
\end{equation}
if $h\geq r$, and in
\begin{equation}
\exp\left(-\sum_{i=h+1}^{r}{\frac{s^i}{i}\left(\sum_{\rho \in \mathcal{Z}(g)}{\frac{1}{\rho^i}}\right)}\right) \prod_{\rho \in \mathcal{Z}(f)\backslash \mathcal{Z}(g)}{\left(1- \frac{s}{\rho}\right)\exp\left({\frac{s}{\rho}+\frac{s^2}{2\rho^2}+...+\frac{s^{h}}{h\rho^{h}}}\right)} \nonumber
\end{equation}
if $r\geq h$. From the definition of the genus, the following two series
\begin{equation}
\sum_{\rho \in \mathcal{Z}(g)}{\frac{1}{\rho^i}}\leq \sum_{\rho \in \mathcal{Z}(f)}{\frac{1}{\rho^i}} \nonumber
\end{equation}
are absolutly convergent for $i \in \left\{r,...,h\right\}$. Hence, in those two cases, the right-hand side is entire and of finite order and this completes the proof. \\
 
We shall now reveal the proof of the main result.
\paragraph{Proof of Theorem \ref{main} :} We suppose that
\begin{equation}
L_1(k-s)=N_1^s\gamma(s)L_1^*(s) \quad (\forall s\in\textbf{C}), \nonumber
\end{equation}
\begin{equation}
L_2(k-s)=N_2^s\gamma(s)L_2^*(s) \quad (\forall s\in\textbf{C}) \nonumber
\end{equation}
and we define
\begin{equation}
\rho(s):=\frac{L_1\left(s+\frac{k}{2}\right)}{L_2\left(s+\frac{k}{2}\right)} \quad \text{and} \quad \rho^*(s):=\left(\frac{N_1}{N_2}\right)^{k/2}\frac{L_1^*\left(s+\frac{k}{2}\right)}{L_2^*\left(s+\frac{k}{2}\right)}. \nonumber
\end{equation}
We will show that $\rho$ belongs to $\mathcal{O}_{N_1/N_2}$.
\begin{enumerate}
\item $\rho$ is a meromorphic function on $\textbf{C}$ (as the quotient of two meromorphic functions) represented by an absolutely convergent ODS.
\item $\forall s\in \textbf{C}$, $\rho(-s)=(N_1/N_2)^s\rho^*(s)$ where $\rho^*$ is also represented by an absolutely convergent ODS.
\item We introduce $f_1(s):=(P_1L_1)(s+k/2)$ (resp. $f_2(s):=(P_2L_2)(s+k/2)$). By assumption, $f_1/f_2$ has finitely many poles : let $\mathcal{Z}$ be the finite multiset of those poles. From Lemma \ref{order} (with $f(s)=f_1(s)$ and $g(s)=f_2(s)/\prod_{\rho \in \mathcal{Z}}{(s-\rho)}$), we know that
\begin{equation}
\prod_{\rho \in \mathcal{Z}}{(s-\rho)}\frac{f_1(s)}{f_2(s)} \nonumber
\end{equation} 
is entire and of finite order in $\textbf{C}$. Thus, $\rho$ is the product of a rational fraction and an entire function since
\begin{equation}
\rho(s)=\left(\frac{P_2(s+k/2)}{P_1(s+k/2)\prod_{\rho \in \mathcal{Z}}{(s-\rho)}}\right)\left(\prod_{\rho \in \mathcal{Z}}{(s-\rho)}\frac{f_1(s)}{f_2(s)}\right). \nonumber
\end{equation}
Taking $M$ to be strictly larger than all the imaginary parts of the poles of the rational fraction, we obtain the third condition that defines $\mathcal{O}_{N_1/N_2}$.  
\end{enumerate}
Hence $\rho\in \mathcal{O}_{N_1/N_2}$. If $N_1/N_2$ is not an integer, $\rho$ is the constant function equal to zero and the result is clear. If $N:=N_1/N_2$ is an integer, there exist complex numbers $(a_d)$ such that
\begin{equation}
\rho(s)=\sum_{u|N}{\frac{a_u}{u^s}} \quad (\forall s\in \textbf{C}) \nonumber
\end{equation}
(by Theorem \ref{even}) and this is the result. \\

Note that the proof infers $N=N_1/N_2$.
\bibliographystyle{plain}
\bibliography{A_Theorem_for_Distinct_Zeros_of_L-Functions}

\begin{thebibliography}{1}

\bibitem{Bom}
A.~Perelli E.~Bombieri.
\newblock Distinct zeros of {L}-functions.
\newblock {\em Acta Arithmetica}, 83(3):271--281, 1998.

\bibitem{Gar}
P.~Garret.
\newblock {W}eierstrass and {H}adamard products.
\newblock {\em Complex analysis course on Weierstrass and Hadamard products},
  pages 1--8, 2015.

\bibitem{Ham}
H.~Hamburger.
\newblock Uber die riemannsche funktionalgleichung der $\zeta$-funktion.
\newblock {\em Mathematische Zeitschrift}, 11(3-4):224--245, 1921.

\bibitem{Rag}
R.~Raghunathan.
\newblock A comparaison of zeros of {L}-functions.
\newblock {\em Mathematical Research Letters}, 6:155--167, 1999.

\end{thebibliography}

\end{document}